\documentclass[10pt]{article}

\usepackage{amsmath}
\usepackage{amssymb}
\usepackage{a4}

\def \dim{\mathop{\rm dim}\nolimits}

\def \mult{\mathop{\rm mult}\nolimits}

\def \halt#1#2 {H^{alt} _{#1} (#2)}
\def \haltz#1#2 {H^{alt} _{#1} (#2;\Z ) }

\def \dk#1 {D^k(#1)}
\def \ep#1 {\varepsilon _{#1}}

\def \constant#1 {{\bf \Q }^\bullet _{#1}}
\def \conshft#1#2 {{\bf \Q }^\bullet _{#1} [#2]}
\def \cplx#1{{\bf #1}^\bullet}

\def \rank{{\mbox{\rm{rank}}}}
\def \sheaf#1{{\underline{#1}}}

\def \R{{\mathbb R}}
\def \C{{\mathbb C}}

\def \Z{{\mathbb Z}}

\def \Q{{\mathbb Q}}
\def \P{{\mathbb P}}

\def \AA{{\cal A}}

\def \LL{{\cal L}}

\def \OO{{\cal O}}

\newtheorem{theorem}{Theorem}[section]
\newtheorem{corollary}[theorem]{Corollary}
\newtheorem{proposition}[theorem]{Proposition}
\newtheorem{lemma}[theorem]{Lemma}
\newtheorem{remark}[theorem]{Remark}
\newtheorem{remarks}[theorem]{Remarks}
\newtheorem{example}[theorem]{Example}
\newtheorem{examples}[theorem]{Examples}
\newtheorem{definition}[theorem]{Definition}

\newenvironment{proof}
 {\begin{trivlist} \item[\hskip \labelsep {\bf Proof.}]}
 {\hfill$\Box$\end{trivlist}}
\begin{document}
\renewcommand{\theenumi}{(\roman{enumi})}
\normalsize

\title{Equisingularity of families of hypersurfaces and applications to mappings} 
\author{Kevin Houston \\
School of Mathematics \\ University of Leeds \\ Leeds, LS2 9JT, U.K. \\
e-mail: k.houston@leeds.ac.uk
}
\date{\today }
\maketitle
\begin{abstract}
In the study of equisingularity of isolated singularities we have the classical theorem of Brian\c{c}on, Speder and Teissier which states that a family
of isolated hypersurface singularities is Whitney equisingular if and only 
if the $\mu ^*$-sequence for a hypersurface is constant in the family.
In this paper we generalize to non-isolated hypersurface singularities. 
By assuming non-contractibility of strata of a Whitney stratification of the non-isolated singularities outside the origin  we show that 
Whitney equisingularity of a family is equivalent to constancy of a certain selection of invariants from two distinct generalizations of the $\mu ^*$-sequence. Applications of this theorem to equisingularity of more general mappings are given.

\end{abstract}

\section{Introduction}
Notions of equisingularity for varieties date back many years, the modern era being started effectively by Zariski. Much work has been done in this area, see \cite{bobadilla} for some recent significant and interesting results. A classical theorem of Teissier, and Brian\c{c}on and Speder gives conditions for equisingularity of a family of complex hypersurfaces such that each member has an isolated singularity.
In this case the family is called {\em{Whitney equisingular}}
if the singular set of the variety formed by the family is a stratum in a
Whitney stratification. 

For any isolated hypersurface singularity we may associate a {\em{$\mu ^*$-sequence}}:
The intersection of the Milnor fibre of the singularity 
and a generic $i$-plane passing 
through the singularity is homotopically equivalent to a wedge of spheres, the number of which is
denoted  $\mu ^i$. This is a sequence of analytic invariants. 

The result of Brian\c{c}on--Speder--Teissier is that 
the $\mu ^*$-sequence is constant in the family if and only if 
the family is Whitney equisingular. 

A natural and long-standing question is what happens in the non-isolated case? 
More precisely, we assume that we can stratify a family so that outside the parameter axis of the family we have a Whitney stratification and seek conditions which give an equivalence between a collection of topological invariants and Whitney equisingularity of the parameter axis. In some sense this was answered in \cite{teissier} using the multiplicity of polar invariants. However, in many situations the number of invariants is very large. We would like a small number of topological or algebraic invariants, defined in a simple manner, which control, and are controlled by, the equisingularity of the family.

An important theorem of Gaffney and Gassler, Theorem 6.3 of \cite{gaffgass}, gives a partial result.
Here they define the sequence $\chi ^*$ as the Euler characteristic of the Milnor fibres that occur for the family. This is an obvious generalization of the $\mu ^*$-sequence as the Euler characteristic of the Milnor fibre is determined  by the Milnor number in the isolated singularity case. The constancy of this sequence does not seem to be sufficient to ensure Whitney equisingularity.
Thus, they define another sequence, the relative polar multiplicities, denoted $m_*$, (see Section \ref{sec:polar_blow} for a precise definition). 
In the case of isolated singularities the constancy of $\mu ^*$ in the family implies the constancy of $m_*$ in the family. 

The Gaffney-Gassler theorem is that Whitney equisingularity of a family implies that the sequences $(m_*, \chi ^*$) are constant in the family. The aim of the paper is to give further conditions to ensure a converse. This is done in Theorem \ref{newthm}. Better than that, the number of invariants is reduced considerably in Theorem \ref{mainthm} to a certain selection from $m_*$ and $\chi ^*$.

The key condition in Theorem \ref{newthm} is that the complex links of strata in the family (outside the parameter axis) have non-trivial homology. There are plenty of examples of spaces with this condition. In fact, in the applications we have in mind there is a plethora of examples. Many more need to be found though.

In Section \ref{sec:setup} we describe the basic notations used and make precise the definition of equisingularity. 
Section \ref{sec:polar_blow} defines the relative polar invariant and Euler characteristic sequences via the method of blowing up of ideals. Two further sequences, the L\^e numbers of Massey, denoted $\lambda ^*$, and Damon's higher multiplicities, denoted $\mu ^*$, are defined in Section \ref{sec:polar_sheaf}. As one can see from the notation, the latter is a generalization of the usual Milnor number. In fact the sequence is equivalent to $m_*$, only the indexing is different. The $\lambda ^*$ sequence is closely related to the $\chi ^*$ sequence - in a family constancy of one implies the constancy of the other. In contrast
to $m_*$ and Damon's $\mu ^*$, however, they are not equal, even after reindexing. Theorem \ref{newthm} gives a partial converse to the Gaffney--Gassler theorem, i.e., gives conditions for $(m_*,\chi ^*)$ constant implies Whitney equisingularity.

The main theorem, Theorem \ref{mainthm}, is given in Section \ref{sec:main}. This gives conditions for equivalence of different sequences and Whitney equisingularity. It also shows that one needs only a selection of invariants from two sequences - we do not require every element from both sequences.

Section \ref{sec:appl} gives applications of the main theorem to families of maps with isolated instabilities such that the discriminant is a hypersurface. In this case we have a large supply of hypersurfaces that satisfy an important condition of the main theorem - the complex link of strata are not homologically trivial. The equisingularity of a family of maps - rather than merely equisingularity of their discriminants - is considered for corank $1$ multi-germs $f:(\C ^n, \underline{x}) \to (\C ^p,0)$. Ultimately, we can produce theorems concerning topological triviality of families of maps.

Note that in the applications we treat the case of multi-germs. 
Despite not requiring much more work and their great importance, particularly in the study of images of maps, these have often been ignored in the past. 

A number of remarks concerning the work of others and areas of possible research are made in the final section.


\section{Equisingularity and basic definitions}
\label{sec:setup}
In this section we give some notation and basic definitions related to equisingularity for the sets
and the complex analytic maps that concern us. 

Standard definitions from Singularity Theory, such as finite 
$\AA $-determinacy, can be found in \cite{dupwall} and \cite{wall}.
The zero set of a map $F$ will be denoted $V(F)$, and its singular set, i.e., the points in the domain where the
rank of the differentiable is less than the codomain, will be denoted by $\Sigma (F)$.
A differentiable map is called {\em{corank 1}} if its differential 
has corank at most 1 at all points. Note that, for convenience, this includes the case of non-singular maps.

Often we shall need to move from a germ and choose a representative, or a
smaller neighbourhood, etc. Since this is entirely standard and is obvious 
when it occurs, no explicit mention shall be made of the details as they 
will be distracting to the exposition.

\begin{definition}
Let $X$ be complex analytic set and $Y$ a subset of $X$. We say that
$X$ is {\em{Whitney equisingular along $Y$}} if $Y$ is a stratum
of some Whitney stratification of $X$. 
\end{definition}
This notion has been the subject of considerable investigation, see \cite{gaffmass}
for a survey from ten years ago and \cite{bobadilla} for more recent developments
in the hypersurface case. In the more general case of maps 
Gaffney made many of the fundamental definitions for the study of the equisingularity,
see \cite{polar}. His work has been continued by him and others, see \cite{gaffpairs,diswe1,perez}.

The famous example of Brian\c{c}on and Speder \cite{bs} shows that even in the
hypersurface case the notion of equisingularity can be a delicate one.

\begin{example}[\cite{bs}]
Let $f(x,y,z,t)=z^5+ty^6z+y^7x+x^{15}$. This is a family of quasihomogeneous 
hypersurface singularities indexed
by $t$ such that $f_t:(\C ^3,0) \to (\C ,0) $ has an isolated singularity at $(0,0,0)$ and the
Milnor number is constant for all $t$.

Since $\Sigma (f)$  is a manifold the obvious stratification of $f^{-1}(0)$ consists of the manifolds $f^{-1}(0) \backslash \Sigma (f)$ and
$\Sigma (f)$.
However, this stratification is not Whitney equisingular along $\Sigma (f)$ as the Whitney
conditions fail at $(0,0,0,0)$.  
What is really interesting is that the family is still topologically trivial.
\end{example}
This example shows that the Milnor number is insufficient to achieve a Whitney stratification; the Brian\c{c}on--Speder--Teissier theorem tells us we need to look at generic slices of the hypersurfaces.
Precise conditions to achieve topological triviality are still the subject of current research.

\section{Polar invariants via blowing up}
\label{sec:polar_blow}
We shall consider what are called {\em{polar invariants}}, these are very important in the study of 
equisingularity, see for example \cite{teissier}. In this section we will consider them as arising from 
the method of blowing up ideals and in the next from the viewpoint of sheaf theory.

Let $f:(\C ^{N+1},0)\to (\C ,0)$ be a complex analytic function, and
denote the Jacobian ideal by $J(f)$:
\[
J(f)=\left( \frac{\partial f}{\partial z_0} , \dots , 
\frac{\partial f}{\partial z_N}  \right)
\]
for coordinates $z_0, \dots , z_N$ in $\C ^{N+1} $.

\begin{definition}
The {\em{blowup of $\C ^{N+1}$ along the Jacobian ideal}}, denoted
$Bl_{J(f)}\C ^{N+1}$, is the closure in $\C ^{N+1} \times \P ^N$ of the
graph of the map
\[
\C ^{N+1} \backslash V(J(f))\to \P ^N , \qquad x\mapsto 
\left( \frac{\partial f}{\partial z_0}(x) : \dots :
\frac{\partial f}{\partial z_N}(x)  \right),
\] 
where $V(J(f))$ is the zero-set of $J(f)$.
\end{definition}

A hyperplane $h$ in $\P ^N$ can be pulled back by the natural
projection $p:\C ^{N+1} \times \P ^N\to \P ^N$ to a Cartier divisor, $H$,  
on $Bl_{J(f)}\C ^{N+1}$, (provided $Bl_{J(f)}\C ^{N+1}$ is not contained in the product of $\C ^{N+1}$
and $h$). 
We call this a {\em{hyperplane on $Bl_{J(f)}\C ^{N+1}$}}.

Let $b:\C ^{N+1} \times \P ^N \to \C ^{N+1}$ be the other natural projection.
For suitably generic hyperplanes $h_1$, \dots , $h_k$ in $\P ^N$, 
the multiplicity at the origin of 
$b(H_1\cap \dots \cap H_k\cap Bl_{J(f)}\C ^{N+1} )$ is well-defined invariant of $f$,
see \cite{gaffgass}.

\begin{definition}
For $1\leq k \leq N$, the {\em{$k$th relative polar multiplicity of $f$}} is the multiplicity
of the scheme $b_*(H_1\cap \dots \cap H_k\cap Bl_{J(f)}\C ^{N+1} )$ at the
origin. It is denoted by $m_k(f)$.
\end{definition}
From this we can define a sequence of invariants $m_*(f)$.
Full details of the above construction and proofs of the various assertions
can be found in \cite{gaffgass} where the authors also show that the situation
can be generalized to ideals other than the Jacobian.

We can now define another, perhaps more familiar, sequence of invariants; these
have a topological nature.
\begin{definition}[\cite{gaffmass} p238]
Let $f:(\C ^{N+1},0) \to (\C ,0) $ be a complex analytic function and
$L^i\subseteq \C ^{N+1}$ be a generic $i$-dimensional linear subspace.  
Denote the reduced Euler characteristic of the Milnor Fibre of $f|L^i$
by $\widetilde{\chi }^i (f)$.
\end{definition}
From this we can define a sequence
\[
\widetilde{\chi }^*(f):= \left( \widetilde{\chi }^{2}(f), \dots , \widetilde{\chi }^{N+1}(f)  \right) .
\]
In the case of an isolated singularity, this (effectively) 
reduces to the standard $\mu ^*$-sequence in Equisingularity Theory.

\begin{remark}
It transpires that the number $\widetilde{\chi }^1 (f)$ is not needed in the theory
in \cite{gaffgass} and so is omitted. This is because in a family of hypersurfaces
$\widetilde{\chi }^2(f)$ of a member will be the Euler characteristic of the Milnor fibre of a
plane curve singularity and so constancy of this implies the constancy of the multiplicity of
the singularity, which implies the constancy of $\widetilde{\chi }^1(f)$.
\end{remark}
In, for example \cite{numcontrol} page 73, Massey shows how one 
can calculate the reduced Euler characteristic
in practice: it is equal to the alternating sum of the L\^e numbers. His definition of L\^{e} numbers
involves taking certain hyperplanes. The precise conditions needed on these hyperplanes 
are not important here, what is important is that they need not be generic. (The lack of
genericity means that we can in practice calculate the L\^e numbers.)

The generic L\^{e} numbers, i.e., those formed by taking generic hyperplanes, 
can be defined using the blowing up setup as follows.
Denote the exceptional divisor of the blowup by $E$. Then, Gaffney and Gassler
in \cite{gaffgass} define the $k$th L\^{e} number, $\lambda _k (f)$, to be the multiplicity of
\[
b_*\left( H_1 \cdots  H_{k-1} \cdot E  \cdot Bl_{J(f)}\C ^{N+1} \right) ,
\]
where $\cdot $ denotes intersection product and $1\leq k \leq N$.
That these coincide with Massey's definition can be found in Theorem II.1.26 on page 61
of \cite{numcontrol}.
Note, however, 
 that this defines the number by codimension, whereas Massey defines the L\^{e} numbers
by dimension. To avoid confusion in this paper for our invariants we will generally use superscripts 
to denote dimension and subscripts to denote codimension.
Hence, in Gaffney-Gassler notation of \cite{gaffgass} $\lambda _i(f)$ is Massey's $\lambda ^{N-i+1}(f)$ (which  we shall define in the next section).

The significance of the invariants $m_*(f)$ and $\widetilde{\chi } ^*(f)$ is made clear in \cite{gaffgass}.
\begin{theorem}[Theorem 6.3 of \cite{gaffgass}]
\label{gaffgass6.3}
Suppose that we have a family of maps $f_t:(\C ^{N+1},0)\to (\C ,0)$. Let 
$F:(\C ^{N+1} \times \C ,0)\to (\C ,0)$ be given by $\overline{F}(x,t)=f_t(x),$ so that
$F(x,t)=(\overline{F}(x,t,),t)$ is a one-parameter unfolding. 

If $V(\overline{F})$ admits a Whitney stratification with 
$T=(\{ 0 \} \times \C ,0) \subset (\C ^{N+1} \times \C ,0)$ as a
stratum, then the map
$t\mapsto (m_*(f_t), \widetilde{\chi }^*(f_t) )$
is constant on $T$.
\end{theorem}
The main aim of this paper is to investigate extra conditions
upon $V(\overline{F})$ which imply that the converse holds. Note that in
\cite{gaffgass} the authors do prove a partial converse in their Theorem 6.2 by showing that the smooth part
of $V(\overline{F})$, the smooth part of its critical locus $\Sigma $ and the components of the singular locus of codimension one in $\Sigma $ are all Whitney regular over the parameter axis.

Note that in the case that $f_t$ is a family of 
isolated hypersurface singularities we have that $\mu ^*(f_t)$ is constant is equivalent to  
$\widetilde{\chi }^*(f_t)$ is constant and these imply that $m_*(f_t)$ is constant. Hence, in this
particular case we know by the Brian\c{c}on-Speder-Teissier result that $(m_*(f_t), \widetilde{\chi }^*(f_t) )$
constant does imply that there is a stratification such that $T$ is a Whitney stratum.

\begin{definition}
Suppose that we have a family of maps, we wish to make precise what it means for an invariant of 
the members to be {\em{constant in the family}}. We shall take it to mean that there is an open contractible
neighbourhood of the origin in the parameter space over which the invariant is constant for
elements of  the family. This definition saves us from constantly referring to the neighbourhood.
\end{definition}


\section{Polar invariants via sheaf theory}
\label{sec:polar_sheaf}
Using intersection theory and sheaf theory Massey has given a different interpretation
of the blowing up we have just seen. The material in this section comes mostly from
\cite{masseyduke, massey, numcontrol}.
The book \cite{numcontrol} in particular contains useful appendices on analytic
cycles, intersection theory, and on vanishing cycles for sheaves.

Suppose that $\cplx{F} $ is a complex of constructible sheaves on an analytic space $X$ and that
$f :X \to \C $ is a complex analytic function. Then we denote the {\em{vanishing cycles
of $\cplx{F} $}} by $\phi _f \cplx{F} $. See \cite{numcontrol} Appendix B for a full definition and
important properties of this complex.

Goresky and MacPherson developed a theory of Morse data on stratified spaces with respect to 
constructible sheaves.
Recall their definition of a non-degenerate conormal vector, see \cite{gm} page 160.
Let $p$ be a point in the stratum $S$ of $X$ and let $T^*_SM$ denote the
set of all covectors $\omega \in T_p^*M$ such that $\omega (T_pS)=0$.

\begin{definition}[\cite{gm} p160]
A plane $Q\subseteq T_p(M)$ is called a {\em{generalized tangent space}} if 
$Q=\lim _{i\to \infty} T_{q_i} S_\alpha $ where $S\subset \overline{S_\alpha }$ and
$q_i$ is a sequence of points in $S_\alpha $ converging to $p$.

Also, the set of {\em{non-degenerate normal covectors}} is the set   
\[
C_S:= \{ \omega \in T_S^*M \, | \, \omega (Q) \neq 0 {\mbox{ for any generalized tangent space }}
Q\neq T_p(S) \} .
\]
\end{definition}

\begin{definition}
Let $X\subset \C ^N$ be a complex analytic space with a Whitney stratification 
$\{ S_\alpha \}$ such that
the strata are connected. Let $\cplx{F}$ be a complex of sheaves 
which is constructible with respect to this stratification.

Let $x$ be a point in the $d$-dimensional stratum $S_\alpha$. 
Let $M$ be a normal slice to $S_\alpha $ at $x$ and $L:(\C ^N,0)\to (\C ,0)$ be a
linear map such that $d_pL$ is a non-degenerate covector. 

Then, the {\em{characteristic normal Morse data}} for the pair $(S_\alpha , \cplx{F})$ is
\[
m(S_\alpha ,\cplx{F})
= (-1)^{N-1}\chi \left( \phi _{L|_X} \cplx{F} \right) _x 
= (-1)^{N-d-1}\chi \left( \phi _{L|_{M\cap X }} \cplx{F}_{|M \cap X } \right) _x ,
\]
where $\chi $ denotes the Euler characteristic of the sheaf (at the point $x\in S_\alpha $).
\end{definition}
In the case that $\cplx{F}$ is the constant sheaf $\cplx{\C }_X $ then
we write $m(S_{\alpha })$ and call it simply the {\em{characteristic normal Morse data of the 
stratum}}. In this case,
\[
m(S_\alpha ) = (-1)^{N-d}
\chi \left( B_\epsilon (x) \cap X \cap M , B_\epsilon (x) \cap X \cap M \cap L ^{-1}(\eta ) \right) ,
\]
where $B_\epsilon (x)$ is a sufficiently small open ball of radius $\epsilon $ centred at $x$, and
$\eta \neq 0$ is also sufficiently small.
Note that Massey uses a different notation, i.e., our $m(S_\alpha ,\cplx{F} )$ is his  
$m_\alpha (\cplx{F} )$,  and our $m(S_{\alpha })$ is his $m_\alpha $.

\begin{definition}
The space $B_\epsilon (x) \cap X \cap M \cap L ^{-1}(\eta )$ in the pair above is called the {\em{complex link of the
stratum $S_\alpha $}}.
\end{definition}

For complete intersections the number $m(S_\alpha )$ is very important. 
\begin{remark}
If $X$ is a complete intersection, then the complex link of a stratum is homotopically equivalent to a wedge of spheres,
see \cite{gm} page 187 or \cite{le1}. Provided $S_\alpha $ is not a `top', non-singular stratum
of $X$ (i.e., a stratum of maximal dimension), 
then $m(S_\alpha )$ is just the number of these spheres. See, for example, Example 6.5
of \cite{masseyduke}.

Thus, in the case of complete intersections $m(S_\alpha )\geq 0$.
\end{remark}

In the general case,
Massey calls strata {\em{visible}} if they have the property that $m(S_\alpha )\neq 0$.
More important
to us are the cases in which this number is positive. This latter property will be a vital assumption in later theorems and their applications and so to 
prevent `empty theorems' we need to produce a significant set of examples for which this holds.

\begin{example}
\label{topstrataeg}
Let $S_\alpha $ be a component of the top strata of $X$. That is, $S_\alpha $ is open in the 
non-singular part of $X$. Then, since the normal slice reduces the normal data to a point, 
the complex link of the stratum is empty and the homology of the normal Morse data is just the
homology of a point. Hence, $m(S_{\alpha })=1$.
\end{example}
We now come to some interesting examples which not only supply plenty of examples, they are useful in applications,
see Section \ref{sec:appl}.
\begin{example}[Theorem 7.3 \cite{clink}]
\label{corank1_clink}
Suppose that $F:(\C ^n,\underline{x} )\to (\C ^p,0) $, $n\geq p$, is a stable, corank 1 map such that $n<p$.
If we stratify the image of $F$ by stable type, then $m(S_{\alpha })=1$
for all strata $S_\alpha $. (Stratification by stable type is described in detail in Section 6 of \cite{polar} and in Section \ref{sec:strat} of this paper.)
\end{example}

\begin{example}[\cite{highmult} p33]
\label{nice_dim_clink}
Suppose that $F: (\C ^n,\underline{x} )\to (\C ^p,0) $ is a stable multi-germ in Mather's nice dimensions (see \cite{nicemather} or \cite{dupwall}).
If we stratify the discriminant of $F$ by stable type, then $m(S_{\alpha })=1$.

This is because of the same reasoning that is behind the previous example, i.e., the complex
link is actually homotopically equivalent to the stabilization of an $\AA _e$-codimension 1 germ. 
\end{example}

\begin{example}
Let $X$ be a hypersurface with an isolated singularity at $x$. For the trivial stratification $\{ X\backslash \{x \} , \{ x \}  \} $
we have $m(\{ x \} ) = \mu (X \cap H) $, where $H$ is a generic
hyperplane. 

\end{example}

We can now state a generalization of the Brian\c{c}on--Speder--Teissier Theorem which is a partial converse to Gaffney and Gassler's theorem.

\begin{theorem}
\label{newthm}
Suppose that $f_t$ is a family as in Theorem \ref{gaffgass6.3}.
Suppose further that $f_t$ is reduced and $X\backslash T$ is Whitney stratified so that the characteristic normal Morse
data $m(S_\alpha )$ is non-zero (and hence positive) for all $S_\alpha \subseteq X\backslash T$.

Then, the following are equivalent.
\begin{enumerate}
\item $\left( m_1(f_t), \dots m_{N}(f_t), \widetilde{\chi }^2(f_t), \dots , \widetilde{\chi }^{N+1} (f_t) \right) $
is constant in the family.
\item The stratum $T$ is Whitney equisingular over all the strata $S_\alpha $. 
\end{enumerate}
\end{theorem}
\begin{proof}
(ii)$\implies $(i): 
This is Theorem 6.3 of \cite{gaffgass}, stated here as Theorem \ref{gaffgass6.3}.

(i)$\implies $(ii):
This follows from the argument of Theorem 6.5 of \cite{gaffmass}:
Consider the stratification of $X\backslash T  $.
The set $\overline{S_{\alpha }}$ is a complex analytic set and thus we can take a 
Whitney stratification $\{ R_\beta \}$ of $X$ such that $\overline{S_{\alpha }}$ is a 
union of strata.

For each $S_\alpha $ there exists a unique $\beta $ such that 
$\overline{S_{\alpha }} = \overline{R_{\beta }}$.
By Chapter 3 of Part III of \cite{numcontrol} the exceptional divisor of the blowup of the
Jacobian ideal of $F$ is, as a cycle, the sum of the projectivisation of the conormal of strata
in the stratification $\{ R_\beta \}$, where each stratum has a multiplicity equal to 
$m(R_\beta )$. Thus, as $m(R_\beta )=m(S _\alpha)\neq 0$, the closure of $R_\beta $ is the
image of a component of the exceptional divisor of the blow-up.

Thus, by Theorem 6.5 of \cite{gaffmass} we know that $R_\beta $ satisfies the Whitney condition
along $T$. Actually, more than this is true as in their proof the authors use Teissier's Theorem 
V.1.2 in \cite{teissier} which states that the non-singular part of the closure of $R_\beta $ is
Whitney over $T$. As $\overline{S_{\alpha }} = \overline{R_{\beta }}$ we deduce that 
$S_\alpha $ is Whitney over $T$. 
\end{proof}

Now recall the definition of {\em{polar varieties}} as described in \cite{masseyduke}.
Let $M$ be the affine space $\C ^{N+1}$, and
let $\underline{z} = (z_0, z_1, \dots , z_N)$ denote a choice of coordinates for $M$.
Define $L^i_{\underline{z}}  :M \to \C ^i $ by $L^i_{\underline{z}}  (z) =(z_0, \dots , z_{i-1})$.

Let $Y$ be an analytic subset of $M$ and $p\in Y$.

\begin{definition}
\label{def:abs-polar}
Suppose that $\dim_\C  {\Sigma } (L^{i+1}_{\underline{z}} |_{Y\backslash \Sigma Y} ) 
\geq i$. Then, the {\em{$i$th absolute polar variety}} with respect to the coordinates $\underline{z}$ 
at the point $p$, denoted $\Gamma ^i_{\underline{z}} (Y)$, is
\[
\Gamma ^i_{\underline{z}} (Y) = {\rm{closure}} \left(
{\Sigma } (L^{i+1}_{\underline{z}} |_{Y\backslash \Sigma Y} ) 
\right)
\]
where $\Sigma (f)$ denotes the critical set of the map $f$ and $\Sigma Y$ denotes the singular set
of the set $Y$. If the dimension condition does not hold, then we define 
$\Gamma ^i_{\underline{z}} (Y) $ to be the empty set.
\end{definition}
If the coordinates are chosen to be sufficiently general, then we get the (generic) absolute
polar varieties of \cite{teissier} and \cite{lete} and we drop the $\underline{z} $ and write
$\Gamma ^i (Y)$.

The following definition arises from Section 7 and Theorem 0.5 of \cite{masseyduke}. The
characteristic polar cycle of a complex of sheaves is defined there in a different way, but is shown 
to be equal to the below in `good' situations. 
\begin{definition}
Suppose that $\cplx{F} $ is a constructible sheaf with respect to the Whitney stratification 
$\{ S_\alpha \} $
of $X\subset M$, where $X$ is a complex analytic set.

The {\em{$k$th characteristic polar cycle of $\cplx{F}$}} (at $p$) is the cycle
\[
\Lambda ^k(\cplx{F})  _p = \sum _{S_\alpha  }  m(S_\alpha ,\cplx{F}) \Gamma ^k (\overline{S_\alpha }) 
\]
where the sum is over all $S_\alpha $ such that $p\in \overline{S_\alpha } $.
\end{definition}

Since the coordinates are generic there is a well-defined multiplicity for 
$\Gamma ^k(\overline{S_\alpha }) $ and so we define the multiplicity of 
$\Lambda ^k(\cplx{F})  $ at $p$, denoted by $\lambda ^k_p(\cplx{F})$ to be
\[
\lambda ^k _p(\cplx{F}) := \mult _p (\Lambda ^k(\cplx{F})_p  ) = \sum _{S_\alpha } m(S_{\alpha } ,\cplx{F} )
\mult _p  \Gamma ^k (\overline{S_\alpha }) .
\]

\begin{example}
\label{constsheaf_eg}
Let $f:(\C ^{N+1},0)\to (\C ,0)$ be a complex analytic map and
let $\cplx{F} $ be the constant sheaf $\cplx{\C } _{V(f)}$ on the hypersurface 
$X=V(f) \subset \C ^{N+1}$.
Then, 
\[
\lambda ^k_0(\C ^\bullet _{V(f)} ) = \sum _{S_\alpha } m(S_\alpha ) \mult _0  \Gamma ^k (\overline{S_\alpha }).
\]
\end{example}
This example is very important as we will see in Lemma \ref{lastmasseyeg} that we can relate these invariants to the relative polar multiplicities of $f$ defined earlier.

We shall drop the reference to $p$ in $\lambda ^k _p(\cplx{F})$ since generally $p$ will be the origin.

\begin{example}[Appendix B \cite{numcontrol}, Chp.\ 10 \cite{massey}]
\label{lenumbereg}
Suppose that $f:(\C ^{N+1}  ,0)\to (\C ,0)$ is a complex analytic map on the manifold $M=\C ^{N+1}$.
Let $\cplx{F}  = \phi_f \cplx{\C } _M $.
Then,
$\lambda ^i(\cplx{F} )$ at $0$ is the $i$th L\^e number of $f$ at $0$,
$\lambda ^i(f)$, as defined by
Massey. (Recall that the Gaffney and Gassler version of this definition stated earlier was 
indexed by codimension
whereas the indexing here is by dimension.)

As the codimension of $J(f)$ in $\C ^{N+1}$ is  at least $2$, we have that $\lambda ^N(f)$ is zero.
This is because the sheaf $\phi_f  \cplx{\C } _M $ is only supported on the 
critical points of $f$.
\end{example}

\begin{remark}
In the example above note that in Appendix B of \cite{numcontrol} and Chapter 10 of \cite{massey} Massey
restricts the sheaf of vanishing cycles to its support and shifts the resulting complex to ensure it is perverse.
\end{remark}

%
%
Recall that if $(X,x)$ is a complete intersection complex analytic set, then the complex link of $x$
is a wedge of spheres of real dimension $\dim _\C X -1$.

Let $f:(\C ^{N+1},0)\to (\C ,0)$ be a complex analytic function.
If $H^i$ is a plane of dimension $i$ through the origin, then $V(f)\cap H^i$ is a complete
intersection.

\begin{definition}[Cf.\ \cite{highmult}]
The {\em{$k$th higher multiplicity}} is the number 
\[
\mu ^k(f) = \dim _\C H _{k-1}(\LL ^k; \C)
\]
where $\LL ^k$ is the complex link of $V(f)\cap H^{k+1}$ at $0$ and $1\leq k \leq N$.
\end{definition}
For sufficiently general $H^k$ this is a well defined invariant of $V(f)$.

These invariants are linked to the relative polar multiplicities by the following lemma.
The three parts of which are effectively from Example 8.4 of \cite{masseyduke}.
\begin{lemma}
\label{lastmasseyeg}
Let $f:(\C ^{N+1},0)\to (\C ,0)$ be a complex analytic function and let $f_t$ be an analytic family of such functions. 
\begin{enumerate}
\item We have $\lambda ^0 (\cplx{\C }_{V(f)} )= \mu ^N (f)$.
\item The numbers $\mu ^{i}(f_t)$ are constant in a family  for all $1\leq i\leq N-r$ if and only if
the numbers
$\lambda ^k(\cplx{\C}_{V(f_t)})$ are constant in a family for all $r+1\leq k \leq N-1$,
where is $r$ is a non-negative integer.

\item For all $1\leq i \leq N$ we have $\mu ^i(f) = m _{N-i+1}(f) $. 
\end{enumerate}
\end{lemma}
\begin{proof}
For parts (i) and (ii) we note Massey's statement in Example 8.4 of \cite{masseyduke}
that 
\begin{eqnarray*}
\lambda ^0(\cplx{\C}_{V(f)}) &=& \mu ^N(f), \\ 
\lambda ^1(\cplx{\C}_{V(f)}) &=& \mu ^{N}(f)+ \mu ^{N-1}(f) \\
\lambda ^2(\cplx{\C}_{V(f)}) &=& \mu ^{N-1}(f)+ \mu ^{N-2}(f) \\
 & & \vdots \\
\lambda ^{i}(\cplx{\C}_{V(f)}) &=& \mu ^{N-i+1}(f)+ \mu ^{N-i}(f) \\
 & & \vdots \\
\lambda ^{N-1}(\cplx{\C}_{V(f)}) &=& \mu ^{2}(f)+ \mu ^{1}(f) \\ 
\lambda ^N(\cplx{\C}_{V(f)}) &=& \mu ^{1}(f)+1.
\end{eqnarray*}

Part (iii) is just the comment from the end of Example 8.4 of \cite{masseyduke}. 
See also  Example 6.5 and 6.10 of the same paper for further information.
\end{proof}

\begin{remark}
Note that (i) and (ii) combine so that 
the numbers $\mu ^{i}(f_t)$ are constant in a family  for all $1\leq i\leq N$ if and only if
the numbers
$\lambda ^k(\cplx{\C}_{V(f_t)})$ are constant in a family for all $0\leq k \leq N-1$.
\end{remark}

\begin{example}
\label{mu_in_isol}
Suppose that $f:(\C ^{N+1},0)\to (\C ,0)$ defines an isolated hypersurface singularity.
Then the lemma shows that $\mu ^k(f)$ coincides with the familiar definition of $\mu ^k(f)$ in the $\mu ^*$-sequence of Teissier
(apart from $\mu ^{N+1}(f)$ which is missing). From this and part (iii) of the lemma we are thus justified in calling Theorem \ref{newthm} a generalization of the 
Brian\c{c}on--Speder--Teissier Theorem.
\end{example}

%
%

\section{Reducing the number of invariants}
\label{sec:main}
We turn once again our attention to the main idea of the paper -- using non-triviality of 
normal Morse data outside of a stratum to give a converse to Gaffney and
Gassler's theorem (stated above as Theorem \ref{gaffgass6.3}). This time we shall add an extra condition to
reduce the number of invariants required in Theorem \ref{newthm} even further.

A further condition required is that outside the stratum of interest the family of maps
is locally trivial over the family's parameter. 
At first sight this may seem a strong condition, but it is found in the main
examples of interest. For example, in the classic Brian\c{c}on--Speder--Teissier result
the family has a line of singularities and outside this line at each point the space is a 
manifold and for each a small neighbourhood
is the product of the parameter axis and a neighbourhood of the point $p$ in the space
above the projection to the axis.

In this section we assume the following.
Let $f:(\C ^{N+1} ,0) \to (\C ,0)$ be a reduced hypersurface and $F(x,t)=(\overline{F} (x,t) ,t)$
be a one-parameter unfolding such that $\overline{F}(x,0)=f(x)$.
Take a representative of $\overline{F}$, also denoted $\overline{F} $, so that
$\overline{F} : U\to \C $ is such that $U\subseteq \C ^{N+1} \times \C $ is an open contractible set.

Let $X=V(\overline{F} )$ and $T= U \cap ( \{ 0 \} \times \C )$, and let 
$\pi : \C ^{N+1} \times \C \to \C $ denote the natural projection. We can identify $T$ with its
image in $\C $ under this map.
Let $\pi ^{-1}(t) = M_t $. For a stratum $S_\alpha $ of a stratification $\{ S_\alpha \}_{\alpha \in \Lambda } $ 
of $X\backslash T$ we define
$S_{\alpha , t} := M_t\cap S_\alpha $. As usual we assume that strata are connected.
\begin{definition}
We say $F$ has a {\em{product structure over $T$}} if the following hold.
\begin{enumerate}
\item The stratification of $X\backslash T$ is Whitney regular with strata of dimension greater than $1$.
\item The set $T$ is contractible.
\item The manifold $M_t$ is transverse to $S_\alpha $ at every point
$(x,t)\in S_\alpha  \subset X \backslash T$ for all strata $S_\alpha $ in the stratification of $X \backslash T$.
\item For all $\alpha \in \Lambda $, $S_{\alpha ,0} \neq \emptyset $.
\end{enumerate}
\end{definition}

Note that 
$M_t=\C ^{N+1} \times \{ t\} $ is a slice such that $\{ S_{\alpha ,t} \} _{\alpha \in \Lambda }$ with 
$\{0 \}$ is a Whitney stratification of $X_t:= M_t \cap X$.

%
%
Due to the product structure on $X$ we can say something about $m(S_\alpha , \cplx{F} )$ for 
various $\cplx {F}$.
\begin{lemma}
\label{productstructure}
Suppose that $F$ has a product structure over $T$. Then 
\begin{eqnarray*}
m(S_{\alpha, 0}) &=& m(S_{\alpha, t})
{\mbox{ and }} \\
m\left( S_{\alpha, 0}, \phi _{f_0} \cplx{\C } _{\C ^{N+1}} \right) &=& 
m\left( S_{\alpha, t}, \phi _{f_t} \cplx{\C }_{\C ^{N+1}} \right)
\end{eqnarray*}
for all $t$ in the family and $S_{\alpha ,t} \neq \{ 0 \} $.
\end{lemma}
\begin{proof}
As $S_{\alpha , t}= S_\alpha \cap M_t$ inherits its stratification from $S_\alpha $ the complex link of
$S_{\alpha ,t }$ is just the complex link of the stratum $S_\alpha $. Therefore,
$m(S_{\alpha ,0}) = m(S_{\alpha ,t })$ as strata are connected and $S_{\alpha ,0} \neq \emptyset $.

For the second part, we note that we have just shown that the characteristic normal Morse data is in effect constant along $T$ (recall that $T$ is contractible) and hence 
we must show that for points $p_t\in S_{\alpha ,t}$ and $p_0 \in S_{\alpha ,0}$ that there exist neighbourhoods
$U_t$ and $U_0$ and a stratum preserving homeomorphism $h:U_t\to U_0$ such that 
\[
h_*\left(  \phi _{f_t} \cplx{\C }_{M_{t}} | _{U_t} \right) \cong 
\phi _{f_0} \cplx{\C }_{M_{0}} | _{U_0} 
\] 
in the bounded constructible derived category. Note that  $\phi _{f_t} \cplx{\C }_{M_{t}}$ is constructible for all $t$ because $\cplx{\C }_{M_{t}}$ is constructible on $M_t$ and $M_t$ can be Whitney stratified so that $V(f_t)$ is a union of strata. Therefore the stratification is Thom $A_{f_t}$, see \cite{bmm} or \cite{parusinski}, and hence by Thom's Second Isotopy Lemma we have the triviality over strata required for constructibility.

Next, the above isomorphism amounts to saying that at every point we can find an 
isomorphism between the vanishing cycles of $f_t$ and $f_0$ at corresponding points in the homeomorphism of
$U_t$ and $U_0$. 

Since $F$ has a product structure over $T$ we have that at every point outside $T$ that the
fibres of $F$ over $T$ are topologically trivial and this homeomorphism is stratum preserving. 
By \cite{le1} we know that topologically equivalent hypersurface singularities have homotopically
equivalent Milnor fibres. Hence, the required result is true. 
\end{proof}
Perhaps the most important fact we can deduce from the assumption of a product structure is that the multiplicity of the absolute polar varieties of the strata of fibres is upper semi-continuous.
\begin{proposition}
\label{prop:semi-cont}
Suppose that $F$ has a product structure over $T$ and the hypersurface defined by $f_t:(\C ^{N+1},0) \to (\C ,0)$ given by $f_t(x)=F(x,t)$ is reduced for all $t\in T$.
Then, for $\overline{S_{\alpha ,t}}\neq \{ 0 \}$ we have that $\mult _0\Gamma ^k (\overline{S_{\alpha ,t}} ) $ is upper semi-continuous for $1\leq k \leq \dim \overline{S_{\alpha ,t}} $.

That is, for sufficiently small $t$ we have 
\[
\mult _0\Gamma ^k (\overline{S_{\alpha ,t}} ) \leq \mult _0\Gamma ^k (\overline{S_{\alpha ,0}} ) 
{\text{ for }} 
1\leq k \leq \dim \overline{S_{\alpha ,0}} .
\] 
\end{proposition}
\begin{proof}
First we need to define the relative $i$th polar variety, denoted $\Gamma ^i(Y,h)$, of a closed complex analytic set $Y\subseteq \C ^{N+1}\times \C $ associated to a complex analytic function $h:Y\to \C $ such that $h^{-1}(t)\subseteq M_t$. It is defined similarly to Definition~\ref{def:abs-polar} and so we shall use the same notation from there. Note that the domain of $L^i$ is still $\C ^{N+1}$.

The set $\Gamma ^i_{\underline{z}} (Y,h)$ is the closure of the union for all $s$ of the set
\[
{\Sigma } \left( L^{i}_{\underline{z}} |_{h^{-1}(s) \backslash \left( h^{-1} (s) \cap  \Sigma Y \right) }  \right) 
\]
provided that the dimension of this closure is greater than or equal to $i$.
If the dimension condition does not hold, then we define 
$\Gamma ^i_{\underline{z}} (Y,h) $ to be the empty set.
As in the remarks following Definition~\ref{def:abs-polar}, by taking sufficiently general projections we can drop the reference to $\underline{z}$ and get 
the relative $i$th polar variety of a closed complex analytic set $Y$ associated to a complex analytic function $h$, denoted $\Gamma ^i(Y,h)$.

In \cite{polar} this set is denoted $P_j(Y,h)$, where $j$ is its codimension in $Y$. Similarly the absolute polar varieties of a set $Z$ is denoted $P_j(Z)$ where $j$ is its codimension in $Z$. 

In our setup we shall have $Y=\overline{S_{\alpha }}$, $h=\pi $ and $Z=\overline{S_{\alpha ,t}}$. 

As we shall require the results from \cite{polar} we explicitly note the connection between the notation there and here. We have 
\[
P_j\left( \overline{S_{\alpha }} ,\pi \right) = 
\Gamma ^{\dim \overline{S_{\alpha }} -j } \left( \overline{S_{\alpha }} ,\pi \right) = 
\Gamma ^{\dim \overline{S_{\alpha ,t}} +1 -j } \left( \overline{S_{\alpha }} ,\pi \right) 
\]
and
\[
P_j\left( \overline{S_{\alpha ,t}}  \right) = 
\Gamma ^{\dim \overline{S_{\alpha ,t}} -j } \left( \overline{S_{\alpha ,t}} \right) . 
\]

In the following we shall assume that $1\leq k \leq \dim \overline{S_{\alpha ,t}} $ without further comment.

By the assumption that $F$ has a product structure we can apply Lemma~5.3 of \cite{polar} and see that 
\[
P_j\left( \overline{S_{\alpha }}, \pi \right) \cap \left( \C ^{N+1} \times \{ t \}  \right) = P_j\left( \overline{S_{\alpha ,t}} \right)  
\]
for all $t$ and for $0\leq j \leq  \dim \overline{S_{\alpha ,t}} -1$. In our notation this is
\[
\Gamma ^{k+1} \left( \overline{S_{\alpha }} , \pi \right) \cap \left( \C ^{N+1} \times \{ t \}  \right) = \Gamma ^k \left( \overline{S_{\alpha ,t}} \right) .
\]
Hence, as the multiplicity of a hyperplane slice of a closed analytic set $Z$ is greater than or equal to the multiplicity of $Z$ by Proposition~7, AC VIII.76 of \cite{bbk}, we have
\[
\mult _0 \Gamma ^{k} \left( \overline{S_{\alpha ,0}} \right) \geq \mult_{(0,0)} \Gamma ^{k+1} \left( \overline{S_{\alpha }} , \pi  \right) .
\]
(On the left hand side we are taking the multiplicity at the origin in $\C ^{N+1}$ and on the right it is the multiplicity at the origin in $\C ^{N+1}\times \C $.)

Next, Teissier's result, Proposition~IV.6.1.1 on page 451 of \cite{teissier}, regarding the upper semi-continuity of the multiplicity of relative polar varieties associated to a map gives
\[
\mult_{(0,t)} \Gamma ^{k+1} \left( \overline{S_{\alpha }} , \pi  \right) \leq 
\mult_{(0,0)} \Gamma ^{k+1} \left( \overline{S_{\alpha }} , \pi  \right)
\] 
for all $t$ in some neighbourhood of $0$ in $\C $.

Finally, for all strata $S_\alpha $ of $X$ we have that $T\subseteq \overline{S_\alpha }$. So, essentially by Proposition~VI.2.1 (p.\ 477) of \cite{teissier}, there exists a contractible open neighbourhood $W$ of $(0,0)$ in $\C ^{N+1}\times \C $ such that $W\cap \left( T\backslash \{ (0,0) \} \right)$ is a Whitney stratum in the obvious stratification of $X\cap W$. (That is, the one given by the stratification of $\left( X\backslash T\right) \cap W $ with the addition of $\left( T\backslash \{ (0,0) \}\right) \cap W $ and  $\{ (0,0) \}$). Hence for all $(0,t)\in T\backslash \{ (0,0) \}$ in this neighbourhood we have by Theorem~5.6 of \cite{polar} that 
\[
\mult _0 \Gamma ^{k} \left( \overline{S_{\alpha ,t}} \right) = \mult_{(0,t)} \Gamma ^{k+1} \left( \overline{S_{\alpha }} , \pi  \right) .
\]
This provides the result
\[
\mult _0\Gamma ^k (\overline{S_{\alpha ,t}} ) \leq \mult _0\Gamma ^k (\overline{S_{\alpha ,0}} ) 
\] 
that we seek.
\end{proof}

Now we can state a lemma which will be at the heart of our next theorem.

\begin{lemma}
\label{mainlemma}
Suppose the following:
\begin{enumerate}
\item The map $F$ has a product structure over $T$.
\item The hypersurface defined by $f_t:(\C ^{N+1},0) \to (\C ,0)$ given by $f_t(x)=F(x,t)$ is reduced.
\item The characteristic normal Morse data $m(S_{\alpha, 0})$ is non-zero (and hence positive) for all $S_{\alpha ,0} \neq \{ 0\}$.
\end{enumerate}
Then, 
\[
\mu ^i (f_t) {\mbox{ is constant for all }} 1 \leq i \leq N-r  
\]
implies that 
\[
\lambda ^i(f_t) {\mbox{ is constant for all }} r+1 \leq i \leq N-1 ,
\]
where $r$ is a non-negative integer.
\end{lemma}
\begin{proof}

From the definition of $\lambda ^k(\cplx{\C }_{V(f_t)})$ (and as the maps $f_t$ are reduced) 
at the origin of $\C ^{N+1} \times \{ 0 \} $ we see that 
\[
\lambda ^k (\cplx{\C }_{V(f_t)} ) = \sum_{S_{\alpha , t}} m(S_{\alpha ,t} ) \mult _0 \Gamma ^k
(\overline{S_{\alpha ,t}} ) . 
\]
Since $\Gamma ^0 (\overline{S_{\alpha ,t}} )=\emptyset $ for all  $\overline{S_{\alpha ,t}} \neq \{ 0 \}$,
and $\Gamma ^k (\{0 \} )=\emptyset $ for all $k \geq 1$, this reduces to 
\begin{eqnarray*}
\lambda ^k (\cplx{\C }_{V(f_t)} ) &=& \sum_{S_{\alpha , t}\neq \{ 0 \} } m(S_{\alpha ,t} ) \mult _0\Gamma ^k (\overline{S_{\alpha ,t}} ) {\mbox{ for }} k\geq 1 , {\mbox{ and }}\\
\lambda ^0 (\cplx{\C }_{V(f_t)} ) &=& m(\{ 0 \}  ) = \mu ^N (f_t) .
\end{eqnarray*}
The last equality above comes from Lemma \ref{lastmasseyeg}(i).

As we have a product structure then
\[
m(S_{\alpha ,t}) = m(S_{\alpha ,0}) 
\]
for all $t$ in the family by Lemma \ref{productstructure}.

As $\mult _0\Gamma ^k (\overline{S_{\alpha ,t}} ) $ is upper semi-continuous by Proposition~\ref{prop:semi-cont}
and $m(S_{\alpha ,t}) \neq 0 $ then we can deduce that for each $k\geq 1$
\begin{eqnarray*}
& &\mult _0 \Gamma ^k (\overline{S_{\alpha ,t}} ) 
 {\mbox{ constant for all }}  \overline{S_{\alpha ,t}} \neq \{ 0 \}     \qquad \qquad \qquad \qquad \qquad \hfill{(*)} \\
& & \iff
\lambda ^k (\sheaf{\C }^\bullet_{V(f_t)} ) {\mbox{ constant.}} 
 \end{eqnarray*}

Now consider the sheaf of vanishing cycles for $f_t$, $\phi _{f_t}  \cplx{\C }_{M_t} $.
Then, for $0\leq k \leq N$, $\lambda ^k (\phi _{f_t}  \cplx{\C }_{M_t} )$ is the L\^e number
of $f_t$, $\lambda ^k(f_t)$, see Example \ref{lenumbereg}, and we have
\[
\lambda ^k ( \phi _{f_t}  \cplx{\C }_{M_t} ) = \sum_{S_{\alpha , t}} m(S_{\alpha ,t} , \phi _{f_t}  \cplx{\C }_{M_t} ) \mult _0 \Gamma ^k
(\overline{S_{\alpha ,t}} ) . 
\]
As we have a product structure we have 
\[
m(S_{\alpha ,t} , \phi _{f_t}  \cplx{\C }_{M_t} )
\]
is constant in the family by Lemma \ref{productstructure}. 

Therefore, for each $k\geq 1$,
\begin{eqnarray*}
& &
\mult _0 \Gamma ^k (\overline{S_{\alpha ,t}} ) 
 {\mbox{ constant for all }}  \overline{S_{\alpha ,t}} \neq \{ 0 \}   \\
& & \implies \lambda ^k ( \phi _{f_t}  \cplx{\C }_{M_t} ) {\mbox{ constant.}}  \qquad \qquad \qquad \qquad \qquad \qquad \qquad \qquad \hfill{(**)}
\end{eqnarray*}
Thus, by  Lemma \ref{lastmasseyeg}(ii), ($*$) and ($**$) we get the statement.
\end{proof}
\begin{remark}
Note that if we have constancy of all the $\mu ^i(f_t)$ for $1\leq i\leq N$, then we control all the L\^{e} numbers 
except $\lambda ^0(f_t)$. In view of the classic Brian\c{c}on--Speder--Teissier result this is not 
surprising. The $\mu ^i(f_t)$ are all the numbers in the classical $\mu ^*$-sequence except the Milnor number of 
the original map $f_t$ and this latter is just $\lambda ^0(f_t)$, see \cite{massey} or \cite{numcontrol}
\end{remark}
\begin{remark}
In the proof of the lemma, a key reason for controlling the higher multiplicities (and hence the 
relative polar multiplicities)
is that the multiplicities of the absolute polar varieties of the
strata are kept constant. It is well known that constancy of these (with some other conditions)
can be used to control the Whitney conditions. See \cite{polar} and \cite{teissier}.
 
This result is perhaps not surprising when one considers one of the main theorems in
\cite{lete}. In Th\'{e}or\`{e}me 4.1.1 of that paper the $\mu ^i (f)$ (and hence 
$\lambda ^i (\cplx{\C }_{V(f)})$)
are connected to the terms $\chi _{d_{\alpha _0}+1} (X, X_{\alpha _0}) $ and 
$\chi _{d_{\alpha _0}+2} (X, X_{\alpha _0}) $, and the $m(S_\alpha )$ correspond to
the $1-\chi _{d_{\alpha } +1 } (X,X_\alpha ) $.

It would be interesting to explore the connection with the work of \cite{lete} and make it more explicit.
\end{remark}

\begin{remark}
\label{rk:vancycpos}
One of the assumptions of the lemma is that the characteristic normal Morse data are positive. This leads to `constancy
of the $\mu ^i$ imply constancy of the $\lambda ^i$'. The same type of proof can be used to show that if the 
$m\left( S_{\alpha, t}, \phi _{f_t} \cplx{\C }_{\C ^{N+1}} \right) $ data are positive, then 
`constancy of the $\lambda ^i$ imply constancy of the $\mu ^i$'.

This may be of interest as cases where the $m\left( S_{\alpha, t}, \phi _{f_t} \cplx{\C }_{\C ^{N+1}} \right) $ 
are positive do occur; for example, the classical Brian\c{c}on--Speder--Teissier result. As there are not many other
obvious examples, we have chosen not to state precisely this version of the lemma. It would, however, be interesting to
find more examples.
\end{remark}

We state another useful lemma for relating invariants in families.
\begin{lemma}
\label{chi_vs_lambda0}
Suppose that $\lambda ^i(f_t)$ is constant for all $1 \leq i \leq N-1$.
Then,
\[
\widetilde{\chi }^{N+1}(f_t) {\mbox{ is constant }}
\iff 
\lambda ^0(f_t) {\mbox{ is constant.}} 
\] 
\end{lemma}
\begin{proof}
By, for example \cite{numcontrol} page 73, the reduced Euler characteristic of the Milnor fibre of $f_t$ (which is
equal to $\widetilde{\chi }^{N+1}(f_t)$) is equal to the alternating sum of the L\^e numbers
$\lambda ^i(f_t)$. From this the lemma follows.
\end{proof}

The main theorem is that we can give a converse to Gaffney and Gassler's Theorem 6.3
(stated above as Theorem \ref{gaffgass6.3}) with fewer invariants. Compare also with Lemma 3.1 in \cite{brazil03}. 
We return to the set up for $f$ and $F$ from the start of this section.
\begin{theorem}
\label{mainthm}
Suppose the following:
\begin{enumerate}
\item The map $F$ has a product structure over $T$.
\item The hypersurface defined by $f_t:(\C ^{N+1},0) \to (\C ,0)$ given by $f_t(x)=\overline{F}(x,t)$ is reduced.
\item The characteristic normal Morse data $m(S_{\alpha, 0})$ is non-zero (and hence positive)
for all $S_{\alpha ,0} \neq \{ 0\}$.
\end{enumerate}
Then, the following are equivalent.
\begin{enumerate}
\item $\left( \mu ^1 (f_t) , \dots , \mu ^{N} (f_t), \lambda ^0(f_t) \right) $ is constant in the family.
\item $\left( \mu ^1 (f_t) , \dots , \mu ^{N} (f_t), \widetilde{\chi }^{N+1}(f_t) \right) $ is constant in the family.
\item $\left( m_1(f_t), \dots m_{N}(f_t), \lambda ^0 (f_t), \dots , \lambda ^{N-1} (f_t) \right) $ is constant in the family.
\item $\left( m_1(f_t), \dots m_{N}(f_t), \widetilde{\chi }^2(f_t), \dots , \widetilde{\chi }^{N+1} (f_t) \right) $
is constant in the family.
\item The stratum $T$ is Whitney equisingular over all the strata $S_\alpha $. 
\end{enumerate}
\end{theorem}
\begin{proof}
(i)$\implies $(iii): 
We have $\mu ^i(f_t)=m_{N-i+1}(f_t)$ for all $1\leq i \leq N-1$ by Lemma \ref{lastmasseyeg}. 
The implication then follows from Lemma \ref{mainlemma} with $r=0$.

(iii)$\implies $(i): This is obvious.

(iii)$\iff $(iv):
This is shown on page 726 of \cite{gaffgass}.

(ii)$\implies $(iii): From Lemma \ref{mainlemma}  we know that $\mu ^i(f_t)$ constant
for all $1\leq i \leq N$ implies  that $\lambda ^k(f_t)$ are constant for $1\leq k \leq N-1$.
From Lemma \ref{chi_vs_lambda0} we deduce that  $\lambda ^0(f_t)$ is constant also.

(iv)$\implies $(ii): This is obvious given $\mu ^i(f_t)=m_{N-i+1}(f_t)$.

(iv) $\iff $ (v): This is Theorem \ref{newthm}.
\end{proof}

\begin{remark}
If $m(S_\alpha )=0$ for the stratum $S_\alpha $, then one can see from the proof of Lemma 
\ref{mainlemma}  that if instead of assuming that $m(S_\alpha )\neq 0$ we assume that $\mult _0 \Gamma ^k (\overline{S_{\alpha ,t}} )$ is constant in the family for
all $k$, then the conclusion of the theorem still holds. 

The statement of such a theorem is obviously ugly and so we choose to omit it. However, 
it is an obvious generalization that may be of some interest in certain cases.
\end{remark}

\begin{remark}
Since $m_i(f_t)=\mu ^{N-i+1}(f_t)$, there exists a number of other obvious equivalences that could
have been stated in the above theorem.
\end{remark}

\begin{remark}
It should be noted that (iii)$\iff $(iv) holds in more generality, see page 726 of \cite{gaffgass}.
\end{remark}

\begin{remark}
In light of Remark \ref{rk:vancycpos}, if we drop the condition that $m(S_{\alpha ,t})>0$ and replace it
with $m\left( S_{\alpha, t}, \phi _{f_t} \cplx{\C }_{\C ^{N+1}} \right) >0$, then we can produce the additional
statement that $\chi ^*(f_t)$ is constant is equivalent to $T$ being a Whitney stratum. This allows us yet 
another way to deduce the classical Brian\c{c}on--Speder--Teissier result and again demonstrates that we should find
more examples where the $m\left( S_{\alpha, t}, \phi _{f_t} \cplx{\C }_{\C ^{N+1}} \right) >0$ condition holds.
\end{remark}

\section{Applications}
\label{sec:appl}
We now apply Theorem \ref{mainthm} to reprove some old results, improve others and to give new ones.
In particular we will consider what happens for equisingularity of families of finitely $\AA $-determined multi-germ maps $f_t:(\C^ n,S)\to (\C ^p,0)$.

\subsection*{The classic Brian\c{c}on-Speder-Teissier result}
The first application is to show that Theorem \ref{mainthm} gives the classic Brian\c{c}on-Speder-Teissier result.
The demonstration of this is included as it is hoped that the proof will shed light on the
application of the theorem to families of finitely $\AA $-determined maps.

\begin{theorem}
Let $F:(\C ^{N+1} \times \C ,0)\to (\C ,0\times 0)$ be family of maps $f_t(x)=F(x,t)$ such that each
$f_t$ defines a reduced isolated singularity at the origin.
Then, the singular set of $F$ is Whitney over the non-singular set if and only if 
$\mu ^*(f_t)$ is constant in the family.
Here $\mu ^i(f_t)$ is the classic Milnor number of $f_t$ restricted to a generic
$i$-plane in $\C ^{N+1}$. 
\end{theorem}
\begin{proof}
As the Milnor number $\mu ^{N+1}(f_t)= \mu (f_t)$ is constant
the set $X= \{ x \in f_t^{-1}(0) {\mbox{ for some }} t \} $ has singular set 
equal to $T= \{ 0 \} \times \C \subset \C ^{N+1} \times \C $.
Thus we can partition $X$ into the manifolds by $\{ X \backslash  T , T \} $.
We have a product structure along $T$ as $X \backslash  T $ is a manifold and
$X\cap ( \C ^{N+1} \times \{ t \} ) = f_t^{-1}(0)$ has a stratification which is obviously 
Whitney.
The normal Morse data of $X \backslash  T $ is equal to 1 by Example \ref{topstrataeg}.

The $\mu ^1 (f_t), \dots  , \mu ^N(f_t)$ of Theorem \ref{mainlemma} are the
usual Milnor numbers by Example \ref{mu_in_isol}. The reduced Euler characteristic
$\widetilde{\chi }^{N+1}(f_t)$ is $(-1)^N\mu ^{N+1}(f_t)$ as the Milnor fibre of $f_t$ is a wedge
of spheres, the number of which is $\mu ^{N+1}(f_t)$. 

Thus by Theorem \ref{mainlemma}, where (ii)$\iff $(v), we deduce the result. 
\end{proof}

\subsection*{Families of finitely $\AA $-determined map germs}
So far the emphasis has been on hypersurfaces. We shall now generalise to a wider class of maps.
Suppose that we have a complex analytic multi-germ $f:(\C ^n,\underline{x} )\to(\C ^p,y)$ where $\underline x = \{ x_1, \dots x_s \}$ is a finite set of points in $\C ^n$. 
Such a map germ is {\em{stable}} at $y$ if all small perturbations of $f$ are $\AA $-equivalent to $f$, i.e., there exist local diffeomorphisms of source and target between the perturbation and $f$. See \cite{dupwall} or \cite{wall} for detailed definitions.
We remark that unfoldings of stable maps are stable.

Let $J(f)$ be the Jacobian of $f$ and let $\Sigma (f) = \{ x\in \C ^n \, | \, \rank \, df < p \} $. This is the {\em{critical set of $f$}}. Define the {\em{discriminant of $f$}}, denoted $\Delta (f)$, to be the image germ of $\Sigma (f)$ under $f$. Note that for $n<p$ this is just the image of $f$.

We say that $f:(\C ^n,\underline{x} )\to(\C ^p,y)$ is {\em{finitely $\AA $-determined at $y$}} if there exists a neigbourhood $U\subseteq \C ^p$ of $y$ such that for all $z\in U\backslash \{ y \} $ the germ 
$f':(\C ^n,f^{-1}(z) \cap \Sigma (f))\to(\C ^p,z)$ is stable. That is, $f$ has an {\em{isolated instability at $y$}}. This definition is analogous to isolated singularity in the case of spaces.

If we have a finitely $\AA $-determined multi-germ $f:(\C^ n,\underline{x})\to (\C ^p,0)$ where $n\geq p-1$, then the
discriminant of $f$ is a hypersurface, see \cite{dupwall} page 446. We apply Theorem \ref{mainthm} in this context, that is, we show when the discriminant of family of maps is Whitney equisingular. Later we will define equisingularity for maps rather than just complex analytic sets.

In order to apply Theorem \ref{mainthm} we stratify the discriminant by stable type. \label{sec:strat}
Section 2.5 of \cite{dupwall} and Section 6 of \cite{polar} are good references for the proofs of the following.

Let $G:(\C ^n,\underline{x} )\to (\C ^p,0)$ be a stable map. 
There exist open sets $U\subseteq \C ^n $ and $W \subseteq \C ^p$ such that 
$G^{-1}(W) = U $ and $G:U\to W$ is a representative of $G$.
We can partition $\Delta (G)$ by stable type. That is, $y_1$ and $y_2$ in $\C ^p$
have the same stable type if $G_1:(\C ^n , G_1^{-1}(y_1)\cap \Sigma(G_1) ) \to (\C ^p,y_1)$ and
$G_2:(\C ^n , G_2^{-1}(y_2)\cap \Sigma(G_2) ) \to (\C ^p,y_2)$ are $\AA $-equivalent. 
These sets are complex analytic manifolds. 

We can take strata in $U\subset \C ^n$ by taking the partition
$G^{-1}(S)\cap \Sigma (G) $, $G^{-1}(S) \backslash  \Sigma (G) )$
and $U \backslash G^{-1}(\Delta (G))$ where $S$ is a stratum in the discriminant.

\begin{definition}[\cite{polar}]
A finitely $\AA $-determined multi-germ $f$ has {\em{discrete stable type}} if there exists a versal unfolding of $f$ in which only a finite number of stable types appear.
\end{definition}
We shall consider two main classes of discrete stable type maps: corank $1$ maps and those in Mather's nice dimensions.
Recall that a map is called {\em{corank $1$}} at a point $x$ if its differential is at most one less than maximal at that point. We say the map is corank $1$ if it is corank $1$ at all points. 
The precise conditions for a map-germ $f:(\C ^n,0)\to (\C ^p,0)$ to be in the nice dimensions are given in \cite{nicemather}. In particular, maps with $p\leq 7$ are in the nice dimensions.

\begin{theorem}
\label{stable_is_canonical}
Suppose that $G:(\C ^n,\underline{x} )\to (\C ^p,0)$ is a stable map and that either of the following holds:
\begin{enumerate} 
\item $G$ is in the nice dimensions, or
\item $G$ is corank $1$ and $n<p$.
\end{enumerate}
Then, the stratification of $G$ by stable type is Whitney regular and any Whitney stratification of 
$G$ is a refinement of this, i.e., this stratification is canonical.
\end{theorem}
\begin{proof}
See Lemma 7.2 of \cite{polar} or Section 2.5 of \cite{dupwall}.
\end{proof}
We shall use this theorem without comment.

\begin{definition}
A stable type is called {\em{$0$-stable}} if the stratification by stable type has a $0$-dimensional stratum.
\end{definition}

\begin{examples}
The Whitney cross-cap $(x,y)\mapsto (x, y^2, xy)$ is $0$-stable. The multi-germ from $(\C ^2 , \{x_1, x_2 ,x_3\} )$
to $(\C ^3,0)$ giving an ordinary triple point is $0$-stable.
\end{examples}
By counting the $0$-stables that appeared in a stable perturbation of a map with an isolated instability
Mond produced interesting and useful invariants of finitely $\AA $-determined maps $f:(\C ^2,0)\to (\C ^3,0)$
in \cite{r2r3}.

Now suppose that we have a multi-germ $f:(\C ^n, \underline{x} )\to (\C ^p,0)$ that has an isolated 
instability at $0\in \C ^p$ and that $F(x,t)=(f_t(x), t)$ is a one-parameter unfolding of $f$ such that $f_0=f$ 
and $f_t(x)=0$ for all $x\in \underline{x}$.

Suppose that we have a representative $F:U\to W$ of $F$ with $F^{-1}(W)=U$. The {\it{parameter axes}} in source and target are
\begin{eqnarray*}
S&:=& (\{ \underline{x} \} \times \C ) \cap U \subset \C ^n \times \C \\
{\mbox{and }}\qquad  T&:=& (\{ 0 \} \times \C )\cap W \subset \C ^p \times \C {\mbox{ respectively}}.
\end{eqnarray*}
We would like stratifications of $F$ so that these are strata. First we can aim to stratify the discriminant of $F$ so that
$T$ is a stratum. If $n\geq p-1$, then the discriminant is a hypersurface and so we can apply Theorem \ref{mainthm}.
A harder problem is to stratify the map itself so that both $S$ and $T$ are strata. 
Here our strategy is to use a stratification of the discriminant and pull it back to one on the source.  
 We could then apply Thom's Second Isotopy Lemma to show that
the family is topologically trivial.

\begin{definition}
We say that the {\em{$0$-stables are constant in the family $f_t$}} if there does not exist a curve $X(t)$ in $\Delta (F)$, the closure of which contains $0\in \C^p$,  such that
$f_t$ has a $0$-stable at $X(t)$.
\end{definition}

\begin{definition}
The {\em{locus of instability}} of $F$ is the set of points 
$(y,t)\in (\C ^p \times \C ,0\times 0)$ such that the map
$F:(\C ^n \times \C , F^{-1}(y,t) \cap \Sigma (F)) \to (\C ^p \times \C , (y,t))$ is not stable
\end{definition}

We can now define the types of unfoldings that we require to apply Theorem \ref{mainthm} to discriminants. This was defined (for mono-germs) by Gaffney in \cite{polar}.  
\begin{definition}
Suppose that $f:(\C ^n , \underline{x} ) \to (\C ^p ,0)$ is finitely $\AA $-determined (i.e., has an isolated instability) and has discrete stable type.
Suppose that $F$ is a one-parameter unfolding with a representative $F:U\to W$ such that $F|\Sigma (F) \cap U \to W $ is proper and finite-to-one, and $F^{-1}(0)\cap \Sigma (F) \cap U = \{ (x_1, 0 ) , \dots , (x_s ,0) \} $. 
 
We call $F$ an {\em{excellent unfolding}} if all the following hold.
\begin{enumerate}
\item $F^{-1}(W)=U$.
\item $F(U \cap \Sigma (F) \backslash S )= W\backslash T $. 
\item The locus of instability is $T$.
\item The $0$-stables are constant in the family.
\item If $n=p$, then the degree of the map $f_t$ is constant in the family.
\end{enumerate}
\end{definition}

\begin{remarks}
\begin{enumerate}
\item Gaffney, in \cite{polar}, calls unfoldings {\em{good}} when all the conditions except the $0$-stable one hold. 
\item These conditions can often be checked by analyzing invariants of the members of the family. See \cite{polar} Proposition 6.6 or Theorem 8.7 for example. See also \cite{excellent} for the case of corank $1$ maps with $n<p$.
\item Representatives $F:U\to W$ such that $F|\Sigma (F) \cap U \to W $ is proper and finite-to-one, and $F^{-1}(0)\cap \Sigma (F) \cap U = \{ (x_1, 0 ) , \dots , (x_s ,0) \} $ can always be found, see \cite{dupwall} page 31.
\end{enumerate}
\end{remarks}

%
%
\subsection*{Main theorem on families of discriminants of map-germs}
We now come to the main theorem in the case that our family of hypersurfaces arises
as the discriminant of the unfolding of a finitely $\AA $-determined map-germ.

We use the following notation.
If $f:(\C ^n,\underline{x})\to (\C ^{p},0)$, $n\geq p-1$, is a finitely $\AA $-determined
multi-germ, then the discriminant is a hypersurface. If $F$ is a one-parameter unfolding of 
$f$ of the form $F(x,t)=(\overline{F}(x,t),t)$, then we shall define $f_t$ to be the family $f_t(x)=\overline{F}(x,t)$ and 
define $g_t:(\C ^p,0)\to (\C ,0)$ to be the 
family of functions defining the discriminants of $f_t$. We can choose $g_0$ reduced, so $g_t$ will be reduced for all $t$ in 
some neighbourhood of $0$.

\begin{theorem}
\label{imagemainthm}
Suppose that $f:(\C ^n,\underline{x})\to (\C ^{p},0)$, $n\geq p-1$, is a finitely $\AA $-determined
multi-germ of discrete stable type and that $F$ is a one-parameter unfolding of 
$f$. Assume that the following hold.
\begin{enumerate}
\item The unfolding is excellent. 
\item The characteristic normal Morse data is non-zero for strata that appear in the stratification by stable types of the discriminant of $F$.
\end{enumerate}
Then, the discriminant of $F$ is Whitney equisingular along the parameter axis $T$ if and only if
the sequence $\left( \mu  ^1(g_t), \dots , \mu  ^{p-1}(g_t), \widetilde{\chi }^p (g_t) \right) $ is constant in the family.
\end{theorem}
\begin{proof}
As $F$ is excellent there are no $1$-dimensional strata other than those contained in the parameter axis.
Furthermore, again since $F$ is excellent, we know from Propositions 6.3, 6.4, and 6.5 of \cite{polar} that 
(i) the stratification of $F|U\backslash F^{-1}(T) \to W\backslash T$ by stable types is a Whitney 
stratification and (ii) the induced stratification of $f_t: U\cap (\C ^n \times \{ t \} ) \to W \cap (\C ^p \times \{ t \} )$
is Whitney and has a product structure over $T$ at the origin. 

Thus, by Theorem \ref{mainthm} applied to the family $g_t$ we get the conclusion. 
\end{proof}

\begin{remarks}
\begin{enumerate}
\item Obviously, by Theorem \ref{mainthm}, other equivalent statements are possible, for example, involving the L\^e numbers of $g_t$. The invariants above were chosen as they are the easiest to define and are clearly topological in nature.
\item Note that in analogy with the Brian\c{c}on--Speder--Teissier result we seem to have the smallest number of invariants possible without making any further assumptions.
\end{enumerate}
\end{remarks}

We can now prove a theorem similar to Theorem 6.6 of \cite{gaffmass}.

\begin{corollary}[Cf.\ \cite{gaffmass}]
Suppose that $f:(\C ^n,\underline{x})\to (\C ^{p},0)$ is a finitely $\AA $-determined
multi-germ  and that $F$ is an excellent one-parameter unfolding of 
$f$. Suppose that $f$ is in Mather's nice dimensions
with $n\geq p$.

Then, the discriminant of $F$ is Whitney equisingular along the parameter axis $T$ if and only if
the sequence $\left( \mu  ^1(g_t), \dots , \mu  ^{p-1}(g_t), \widetilde{\chi }^p (g_t) \right) $ is constant in the family.
\end{corollary}
\begin{proof}
Since $f$ is in the nice dimensions it is of discrete stable type.
The stable types appearing in any unfolding will obviously be in the nice dimensions and so 
the complex links of the stable types are non-contractible by Example \ref{nice_dim_clink}. Thus the characteristic normal data is non-zero and so condition (ii) of the theorem is satisfied.
\end{proof}

We can now state a theorem for the case of images with $p=n+1$.
\begin{corollary}
\label{corank1_cor}
Suppose that $f:(\C ^n,\underline{x})\to (\C ^{n+1},0)$ is a corank $1$ finitely $\AA $-determined
multi-germ and that $F$ is an excellent one-parameter unfolding of 
$f$. 

Then, the image of $F$ is Whitney equisingular along the parameter axis $T$ if and only if
the sequence $\left( \mu  ^1(g_t), \dots , \mu  ^{n}(g_t), \widetilde{\chi }^{n+1} (g_t) \right) $ is constant in the family.
\end{corollary}
\begin{proof}
The map $f$ is of discrete stable type because $f_0$ is corank $1$.
Furthermore, the stable types appearing in any unfolding will also be corank $1$.
The complex links of the stable types are non-trivial by Example \ref{corank1_clink} and so condition (ii) of the theorem is satisfied.
\end{proof}

\begin{remark}
See \cite{excellent} for conditions on members of the family to show that $F$ is an excellent unfolding. 
\end{remark}

\subsection*{Main theorem on families of map-germs}
We can stratify the map $F$ so that the parameter axes $S$ and $T$ are strata.
Gaffney initiated this study of equisingularity of finitely $\AA $-determined maps - rather than just hypersurfaces - in \cite{polar}.
His statements were for mono-germs but the extension to multi-germs is fairly straightforward. 
\begin{definition}
Let $F:(\C ^n\times \C ,\underline{x} \times 0)\to (\C ^p\times \C ,0\times 0)$
be a family of maps $F(x,t)=(f_t(x),t)$ such that each 
$f_t:(\C ^n , \underline{x})\to (\C ^p,0)$ has an isolated instability
at the origin.

We say that $F$ is {\em{Whitney equisingular (along the parameter axes)}} if there is a representative
$F:U\to W$ so that
$U\subseteq \C ^n\times \C $ and $W\subseteq \C ^p\times \C$ can be Whitney stratified so that
\begin{enumerate}
\item $F$ satisfies Thom's $A_F$ condition, and
\item the parameter axes $S=\{ \underline{x} \} \times \C  \subseteq \C ^n \times \C $, and
$T=\{ 0\} \times \C  \subseteq \C ^p \times \C $ are strata. That is, the parameter axes are strata. 
\end{enumerate}
\end{definition}


\begin{remark}
By Thom's Second Isotopy Lemma if a family is Whitney equisingular, then it is topologically trivial.
\end{remark}
For mono-germs we can improve on the main theorem in \cite{brazil03}.
\begin{theorem}[Cf.\ \cite{brazil03} Theorem 3.3]
Suppose that $f:(\C ^n,0)\to (\C ^{n+1},0)$ is a corank $1$ finitely $\AA $-determined
mono-germ and that $F$ is an excellent one-parameter unfolding of 
$f$. 

Then, $F$ is Whitney equisingular along the parameter axes $S$ and $T$ if and only if
the sequence $\left( \mu  ^1(g_t), \dots , \mu  ^{p-1}(g_t), \widetilde{\chi }^p (g_t) \right) $ is constant in the family.
\end{theorem}
\begin{proof}
Clearly, if $F$ is Whitney equisingular along $S$ and $T$, then the image is Whitney equisingular along $T$ and hence by
Corollary \ref{corank1_cor} the sequence is constant in the family.

For the converse the same corollary implies that if the sequence is constant, then the image is Whitney equisingular 
along $T$. The main theorem of \cite{gaffcorank1} implies that the source may also be Whitney stratified so that $S$ is a stratum. Since Gaffney only proves this for mono-germs we also have to restrict to mono-germs.

The Thom $A_F$ condition follows automatically because if $h:Y \to Z$ is a finite complex analytic map with $Y$ and $Z$ Whitney stratified so that strata map to strata by local diffeomorphisms, then $h$ satisfies the Thom $A_f$ condition as the kernels in the definition of Thom $A_F$ are all $\{ 0\}$. The map $F$ is finite and so the submersions formed by taking restrictions to strata are in fact local diffeomorphisms.

Hence, $F$ is Whitney equisingular.
\end{proof}

\section{Final remarks}
\begin{remark}
In \cite{polar}, in particular Proposition 8.4, 8.5 and 8.6, there are a number of formulae relating various polar multiplicities and other invariants. The polar multiplicities appear with coefficient equal to $\pm 1$ in all these propositions. This same behaviour
can be seen in the work of Jorge P\'erez, \cite{perez,perez2} and Saia \cite{jps}.

That these coefficients are equal to $\pm 1$ appear to be a reflection of the fact that in the case of $\C ^n $ to $\C ^p$ with $n<p$ and corank $1$ the stable types appearing have characteristic normal data equal to $1$ by Example \ref{corank1_clink}. 

More specifically, the alternating sum of the multiplicities of the characteristic polar cycle of the constant sheaf $\cplx{\C } _{V(g_t)}$ in Example \ref{constsheaf_eg} can often be written as some other well-known invariant, for example a Milnor number. As 
\[
\lambda ^k(\C ^\bullet _{V(g_t)} ) = \sum _{S_{\alpha ,t}} m(S_{\alpha ,t}) \mult _0  \Gamma ^k (\overline{S_{\alpha ,t} })
\]
and since $m(S_{\alpha ,t} )=1$ (in the notation of Lemma \ref{mainlemma}) we get an alternating sum of the polar multiplicities
$\mult _0  \Gamma ^k (\overline{S_{\alpha ,t} })$.

Alternatively, it is well known that the alternating sum of polar multiplicities is equal to the Euler obstruction, see \cite{lete2}. Thus corank 2 maps, and in particular their characteristic normal data, will need to be studied to determine the precise explanation.
\end{remark}

\begin{remark}
As remarked earlier, the analogy with the Brian\c{c}on--Speder--Teissier Theorem shows that the number of invariants in Theorem \ref{mainthm} cannot be reduced any further without extra conditions being imposed. More than this, it seems likely that one needs the non-triviality of the complex links in the theorem. If one has a contractible complex link of some stratum, then it is probable that one can create examples where the $(\mu ^*, \widetilde{\chi }^{N+1})$ sequence is constant but the parameter axis does not satisfy the Whitney conditions. Such an example therefore needs to be found.

It is difficult to express succinctly the reasons behind this strong probability in the current space and so the interested reader is directed to \cite{numcontrol} Part III section 4 and in particular Proposition 4.10.
\end{remark}

\begin{remark}
While it is satisfying to reduce the number of required invariants to $p$ it is unsatisfactory that they are not defined consistently.  In 
$\left( \mu  ^*(g_t), \widetilde{\chi }^p (g_t) \right) $ we have a mix of higher multiplicities and an
Euler characteristic. On the other hand, many other theorems give an even less consistent mix of polar multiplicities, Milnor numbers, L\^e numbers and so on.

However, for equisingularity of maps it is possible to define yet another sequence of invariants using the {\em{disentanglement}} of a map, see \cite{diswe1} and \cite{highmult} Section 4 for a discussion. This sequence is denoted by $\mu _I ^i(f_t)$ where $1\leq i \leq p$ as it depends on the map $f_t$ and not the function $g_t$ defining the discriminant. It is possible to show that $\mu ^i (g_t) =\mu ^i_I(f_t)$ for $1\leq i \leq p-1$. In low dimensional examples it is possible to see that $\widetilde{\chi }^{p+1}(g_t)$ and $\mu _I^p(f_t)$ are connected though not equal, see \cite{diswe1}. Thus Whitney equisingularity of these maps is controlled by the $p$ invariants $\mu _I^*(f_t)$. It would be interesting to prove in general that 
 constancy of $(\mu ^*(g_t), \widetilde{\chi }^{p}(g_t))$ is equivalent to the constancy of $\mu ^*_I(f_t)$.

Another reason for studying this is that $\mu _I^p(f_t)$ is involved in the control over an unfolding being excellent, see \cite{diswe1} and \cite{excellent}.
\end{remark}

\begin{remark}
Also of interest is to find when equisingularity of the discriminant implies equisingularity of the map. Gaffney's theorem of \cite{gaffcorank1} shows that for a family of corank $1$ maps $f_t:(\C ^n,0)\to (\C^{p},0)$, $n<p$, that if the image is Whitney equisingular, then the family is Whitney equisingular. It would be good to know how general this is. 
It seems unlikely, particularly for $n<p$, that an unfolding should have a source that is not Whitney stratified with $S$ a stratum such that the image is Whitney stratified with $T$ a stratum. That is, one would expect an image to be more complicated than its source and the map should not `repair' faults with stratifications.

Furthermore, as Gaffney pointed out to me, in the case of maps, one cares about the topological triviality and so one need not be restricted to Whitney stratifications. One could use the $c$-regular stratifications of Bekka (see \cite{BekkaWark}) for source and target as this would imply topological triviality of the family by Thom's Second Isotopy Lemma (since the lemma hold for these types of stratifications). Alternatively, one could attempt to find conditions so that Whitney equisingularity of the discriminant implies $c$-regularity of the source (since Whitney stratification is stronger than $c$-stratification, this would again imply topological triviality).
\end{remark}

\end{document}